\documentclass[smallextended,envcountsect]{svjour3}
\usepackage{graphicx,color,cite}
\usepackage{mathrsfs}
\usepackage{amsmath, amscd, amsfonts, amssymb, graphicx, color}
\usepackage[bookmarksnumbered, colorlinks, plainpages]{hyperref}
\usepackage{mathptmx}   
\usepackage[colorinlistoftodos]{todonotes}
\usepackage[left=4cm, right=4cm, top=3.5cm, bottom=3.5cm]{geometry}

\usepackage{latexsym}

\if{
\newtheorem{theorem}{Theorem}[section]

\newtheorem{lemma}{Lemma}[section]

\newtheorem{corollary}{Corollary}[section]
\newtheorem{proposition}{Proposition}[section]

}\fi

\newcommand{\N}{\mathbb N}

\renewcommand\theenumi{\roman{enumi}}

\newcounter{mycount}

\smartqed

\usepackage{etex}

\usepackage{mathtools}

\usepackage{enumitem}
\renewcommand\theenumi{(\roman{enumi})}
\renewcommand\theenumii{(\alph{enumii})}

\usepackage{csquotes}

\usepackage{todonotes}

\makeatletter
\let\orgdescriptionlabel\descriptionlabel
\renewcommand*{\descriptionlabel}[1]{
 \let\orglabel\label
 \let\label\@gobble
 \phantomsection
 \edef\@currentlabel{#1}
 \let\label\orglabel
 \orgdescriptionlabel{#1}
}
\def\th@plain{
 \thm@notefont{}
 \itshape
}
\def\th@definition{
 \thm@notefont{}
 \normalfont
}

\g@addto@macro\th@definition{\thm@headpunct{}}
\g@addto@macro\th@plain{\thm@headpunct{}}
\makeatother

\usepackage{subfig}

\usepackage[final]{showkeys}

\usepackage{etoolbox}

\usepackage{mathrsfs}

\usepackage[titletoc]{appendix}
\usepackage[doc]{optional}

\usepackage{soul}

\usepackage{colortbl,booktabs,
multirow}
\usepackage{xcolor}

\usepackage{cancel}

\usepackage{empheq}

\definecolor{myblue}{rgb}{.8, .8, 1}
\newcommand*\mybluebox[1]{
\colorbox{myblue}{\hspace{1em}#1\hspace{1em}}}

\usepackage{hyperref}
\hypersetup{
colorlinks=true,
linkcolor=blue,
citecolor=blue,
filecolor=magenta,
urlcolor=cyan
}

\usepackage[
open,
openlevel=2,
atend,
numbered
]{bookmark}

\usepackage[capitalize, nameinlink, noabbrev]{cleveref}
\crefname{equation}{}{}
\crefname{chapter}{Chapter}{Chapters}
\crefname{item}{item}{items}
\crefname{figure}{Figure}{Figures}
\crefname{theorem}{Theorem}{Theorems}
\crefname{lemma}{Lemma}{Lemmas}
\crefname{proposition}{Proposition}{Propositions}
\crefname{corollary}{Corollary}{Corollarys}
\crefname{definition}{Definition}{Definitions}
\crefname{fact}{Fact}{Facts}
\crefname{example}{Example}{Examples}
\crefname{algorithm}{Algorithm}{Algorithms}
\crefname{remark}{Remark}{Remarks}
\crefname{note}{Note}{Notes}
\crefname{notation}{Notation}{Notations}
\crefname{case}{Case}{Cases}
\crefname{exercise}{Exercise}{Exercises}
\crefname{question}{Question}{Questions}
\crefname{claim}{Claim}{Claims}
\crefname{enumi}{}{}

\usepackage{pgf}

\usepackage{array}
\usepackage{tabu}

\allowdisplaybreaks

\numberwithin{equation}{section}

\renewcommand\theenumi{(\roman{enumi})}
\renewcommand\theenumii{(\alph{enumii})}

\spnewtheorem*{Proof}{Proof.}{\bf}{\rm}

\begin{document}

\title{Primal Characterizations of Error Bounds for Composite-convex Inequalities\thanks{Research  of the first author was supported by  the National Natural Science Foundations of  China (Grant Nos. 11971422 and 12171419), and funded by Science and Technology Project of Hebei Education Department (No. ZD2022037) and the Natural Science Foundation of Hebei Province (A2022201002). Research of the second author benefited from the support of the FMJH Program PGMO and from the support of EDF.}}

\titlerunning{Characterizations of Error Bounds of Composite-convex Inequalities}

\author{Zhou Wei  \and Michel Th\'era \and Jen-Chih Yao}

\institute{Zhou Wei\at Hebei Key Laboratory of Machine Learning and Computational Intelligence \& College of Mathematics and Information Science, Hebei University, Baoding, 071002, China\\ \email{weizhou@hbu.edu.cn}\\
Michel Th\'era \at XLIM UMR-CNRS 7252, Universit\'e de Limoges, Limoges, France  \\ and \at Federation University Australia, Ballarat\\\email{michel.thera@unilim.fr}\\
Jen-Chih Yao \at Department of Applied Mathematics, National Sun Yat-sen University, Kaohsiung,
Taiwan\\ \email{yaojc@math.nsysu.edu.tw}
}

\date{Received: date / Accepted: date}
\dedication{Dedicated to Roger J-B Wets on the occasion of his 85th birthday. Roger's pioneering work helped stochastic optimization become the stronghold it is today.}

\maketitle

\begin{abstract}
This paper is devoted to primal conditions of error bounds for a general function. In terms of Bouligand tangent cones, lower Hadamard directional derivatives and the Hausdorff-Pompeiu excess  of subsets, we provide several necessary and/or sufficient conditions of error bounds with mild assumptions. Then we use these primal results to characterize error bounds for composite-convex functions (i.e. the composition  of a convex function with a continuously differentiable mapping). It is proved that the primal characterization of error bounds can be established via Bouligand tangent cones, directional derivatives and the Hausdorff-Pompeiu excess if the mapping is metrically regular at the given point. The accurate estimate on the error bound modulus is also obtained.



\end{abstract}

\section{Introduction}
The main goal of this paper is to study error bounds for an inequality defined by a composite-convex function; i.e. the  composition  of a convex function with a continuously differentiable mapping. Error bounds of extended-real-valued functions have been intensively studied for more than half a century. The starting point of the theory of error bounds goes back to the fundamental works by Hoffman \cite{MP31} and Lojasiewicz \cite{Lojasiewicz}. Their results were extensively studied by many authors (cf. \cite{AC1988,Jour2000,KL1999,LP1997,Pang1997,Penot,Rob73, abassi-thera1,abassi-thera2}) and there have been significant developments on error bounds for convex and non-convex functions in recent years. The readers are invited to consult bibliographies \cite{AC,BK,BD1,CK2020,CM2008,FHKO2010,ioffe-JAMS-2,KLT2018,KNT2010,TsL92,MP52,NT2004,NT2008,NT2009,WZ2021,Za} and references therein for theory and applications of error bounds for more details.

Error bounds have played an important role in various aspects  of optimization and variational analysis  including for instance,  sensitivity analysis of linear programming (cf. \cite{Rob73,Rob77}), convergence analysis of descent methods (cf. \cite{Gul92,HLu,IuD90,TsL92,LT,TsB93}), the so-called feasibility problems (cf. \cite{BB,BT,BK}), the domain of image reconstruction (cf. \cite{Combettes}) and many others. Error bounds are closely related with other notions well  known  and used in convex analysis and approximation theory such as  the basic constraint qualification,  the strong conical hull intersection property, the linear regularity, the Abadie constraint qualification and optimality conditions (cf. \cite{DR2,LP1997,Li,MP52}).  Also error bounds are extensively discussed with weak sharp minima of functions, metric subregularity as well as calmness of multifunctions (cf. \cite{AC,BD1,BD2,Gf,ioffe-JAMS-1,ioffe-JAMS-2,M1} and references therein for more details).

When dealing with error bounds, a large   literature is devoted to provide dual characterizations and criteria in terms of subdifferentials or normal cones. To the best of our knowledge, \cite{Io} was among  of the first papers of such kind to state sufficient conditions for error bounds of a constraint system in terms of  the Clarke subdifferential. In 1997 Lewis and Pang \cite{LP1997} studied error bounds for convex inequality systems and provided necessary conditions via subdifferentials and normal cones. In 2003 Ngai and Th\'era \cite{NT2004} provided an error bound estimate and an implicit multifunction theorem in terms of smooth subdifferentials and abstract subdifferentials. In 2004 Zheng and Ng \cite{ZN2004} proved dual characterizations of error bounds for convex inequalities in terms of subdifferentials and normal cones. In 2010 subdifferential characterizations of stability of error bounds for convex constraint inequalities were given in \cite{KNT2010,NKT2010}. In 2018 Kruger, L\'opez and Th\'era \cite{KLT2018} extended results in \cite{KNT2010,NKT2010} and provided subdifferential characterizations of stability of error bounds for convex inequalities in the Banach space setting. It is noted that a pretty natural idea is to study error bounds in terms of various primal derivative-like objects such as directional derivatives, contingent cones or slopes. Several criteria for error bounds were worked out in \cite{CK2020,CM2008,FHKO2010,NT2008,NT2009} based on the primal-type estimate. In terms of contingent cones and directional derivatives, Wei, Yao and Zheng \cite{WZ2014} proved primal characterizations of error bounds for a convex inequality (see \cite[Proposition 5.3]{WZ2014}). Recently, the authors \cite{WZ2021} further studied error bounds of the convex inequality in terms of the Hausdorff-Pompeiu excess  (of subsets), Bouligand tangent cones and directional derivatives, and provided the accurate primal estimate on the error bound modulus (see \cite[Theorem 5.1]{WZ2021}). Based on the works in \cite{WZ2014,WZ2021},  a natural issue is to extend primal results on characterizations of error bounds and on the error bound  modulus by dropping the convexity assumption. Inspired by this issue, our goal in this article is to discuss error bounds of the inequality defined by a composite-convex function. Our work is to provide primal characterizations of error bounds and the accurate estimate on the error bound modulus in terms of  the  notions of Hausdorff-Pompeiu excess  of a set beyond another one,  of the Bouligand tangent cone and of the lower  Hadamard directional derivative.

The paper is organized as follows. In Section 2, we give some definitions and preliminary results. In Section 3, we consider two concepts of metric regularity (of multifunctions) and the Shapiro first order contact property  that are used in our analysis. Section 4 is devoted to the study of error bounds for the inequality defined by a composite-convex function. We first consider error bounds for a general inequality defined by a proper lower semicontinuous (not necessarily convex) function with the Shapiro first order contact property, and provide sufficient and/or necessary primal conditions of error bounds in terms of the Bouligand tangent cone, the  lower  Hadamard directional derivative and the Hausdorff-Pompeiu excess of subsets (see Theorems 4.1 and 4.2). When these results are applied to   error bounds   of  composite-convex inequalities, the primal results on error bounds and  on the error bound modulus can be obtained (see Theorems 4.4 and 4.5).





\section{Preliminaries}

Let $\mathbb{X}$ be a Banach space (or Euclidean space). Let $\mathbf{B}_{\mathbb{X}}$ denote the closed unit ball of $X$. For $\bar x\in \mathbb{X}$ and $\delta>0$, let
$\mathbf{B}(\bar x, \delta)$ denote the open ball with center $\bar x$ and
radius $\delta$. For a subset $\Omega$ of $\mathbb{X}$, we denote by ${\rm cl}(\Omega)$, ${\rm int}(\Omega)$ and ${\rm bd}(\Omega)$ the closure, the interior and the boundary of $\Omega$, respectively.

Let $A$ be a closed subset of  $\mathbb{X}$ and $\bar x\in A$. We denote by
$$
\mathbf{T}^\mathbf{B}(A, \bar x):=\mathop{\rm Limsup}\limits_{t\rightarrow 0^+}\frac{A-\bar x}{t}
$$
the {\it Bouligand tangent cone} (also called  \textit{contingent  cone}) of $A$ at $\bar x$. Thus, $v\in \mathbf{T}^\mathbf{B}(A, \bar x)$ if and only if there exist a sequence $\{v_n\}$ in $X$ converging to $v$ and a sequence
$\{t_n\}$ in $(0,\;+\infty)$ decreasing to 0 such that $\bar
x+t_nv_n\in A$ for all $n\in\mathbb{N}$, where $\mathbb{N}$ denotes
the set of all natural numbers.

For any subsets $C$ and $D$ of $\mathbb{X}$, the {\it excess} of $C$ beyond $D$ is defined as:
\begin{equation}\label{2.2}
  \mathbf{e}(C,D):=\sup_{x\in C} \mathbf{d}(x,D),
\end{equation}
where $\mathbf{d}(x,D):=\inf\{\|x-y\|:y\in D\}$ and the convention is used that
\begin{equation}\label{2.3}
  \mathbf{e}(\emptyset, D):=\left\{
\begin{array}{cl}
 0, & {\rm if} \ D\not=\emptyset,\\
\infty, & {\rm otherwise}.
\end{array}
\right.
\end{equation}


The following proposition provides a characterization for the Hausdorff-Pompeiu excess. We refer the reader to  \cite[Page 138]{DR2}  (or \cite[Proposition 3.1]{WZ2022}) for the proof in details.

\begin{proposition}\label{proposition:pro1}
Let $C$ and $D$ be subsets of $\mathbb{X}$. Then
\begin{equation}
  \mathbf{e}(C,D)=\inf\{\tau\geq 0: C\subseteq D+\tau \mathbf{B}_\mathbb{X}\}.
\end{equation}
\end{proposition}

Given an extended-real-valued   lower semicontinuous function $\varphi:\mathbb{X}\rightarrow \mathbb{R}\cup
\{+\infty\}$, we denote by ${\rm dom}(\varphi):=\{u\in \mathbb{X}:\varphi(u)<+\infty\}$ its domain and by $$
{\rm epi}(\varphi):=\{(x,\alpha)\in \mathbb{X}\times \mathbb{R}:\varphi(x)\leq\alpha\}
$$ its epigraph. $\varphi $ is said to be proper if its domain is nonempty.

The {\it  lower  Hadamard directional derivative} of $\varphi$ at $x$ along direction $h$ is defined as
\begin{equation}\label{2.1}
 \varphi_H^{\prime}(x;h):=\liminf_{t\rightarrow
0^+, h'\rightarrow h}\frac{\varphi(x+th')-\varphi(x)}{t}.
\end{equation}

We close this section with the following lemma cited from \cite[Theorem 4.1]{SYZ}.

\begin{lemma}\label{lemma:lema}
Let $\mathbb{X}$ be a Banach space and  $\Omega$ be a nonempty closed
subset of $\mathbb{X}$. Let $\gamma\in (0,\;1)$. Then for any $x\not\in \Omega$
there exists $z\in \Omega$ such that
\begin{empheq}[box =\mybluebox]{equation*}\gamma\|x-z\|<\min\{\mathbf{d}(x, \Omega), \mathbf{d}(x-z,\mathbf{T}^\mathbf{B}(\Omega,z))\}.\end{empheq}
In particular, if $\mathbb{X}$ is of finite dimension, $z$ can be chosen as the projection of $x$ onto $\Omega$.
\end{lemma}

\setcounter{equation}{0}
\section{Metric regularity and  the Shapiro first order contact property}

This section recalls  the two important concepts of {\it metric regularity} and of  {\it Shapiro first order contact property} that are used in our analysis. It is known that metric regularity of multifunctions occurs to be closely related to Lipschtizian properties of inverse mappings. This well-known and significant property is an extension of surjectivity to nonlinear/set-valued mappings and goes back to the Banach-Schauder open mapping theorem and  to the Lyusternik-Graves theorem. We first recall the definition of metric regularity.\\[1pt]

\noindent{\bf Definition 3.1.} {\it {\rm(i)} Let $F:\mathbb{X}\rightrightarrows \mathbb{Y}$ be a multifunction between two Banach spaces and $\bar x\in \mathbb{X}$. Recall that $F$ is said to be metrically regular at $\bar x$ for $\bar y\in F(\bar x)$, if there exists a constant $\kappa\in (0,+\infty)$ along with neighborhoods $U$ of $\bar x$ and $V$ of $\bar y$ such that}
\begin{empheq}[box =\mybluebox]{equation}\label{3-1}
  \mathbf{d}(x, F^{-1}(y))\leq\kappa \mathbf{d}(y, F(x))\ \ {\it for\ all} \ (x,y)\in U\times V.
\end{empheq}
(ii) {\it A single-value mapping $\varPsi :\mathbb{X}\rightarrow \mathbb{Y}$ is said to be metrically regular at $\bar x$ if $\varPsi $ is metrically regular at $\bar x$ for $\varPsi (\bar x)$.}\\[1pt]

It is known from \cite{Asen-book, M2, Penot2013, ioffe-book, thibault} that metric regularity of a multifunction is proved to be equivalent to the covering property. We refer for instance  the readers to \cite[Theorem 1.52]{M2} for more details on this equivalence relationships and the modulus estimates for metric regularity and the covering property. Further the readers are invited to consult \cite{ioffe-JAMS-1,ioffe-JAMS-2} for a survey and the development on metric regularity and  \cite{NTVT} for some regular properties of graphical tangent and normal cones to paraconvex multifunctions.

From \cite[Lemma 1.56]{M2},  we always have  the automatic closedness of the derivative image for  single-valued  metrically regular  mappings. The following proposition, as a result of independent interest, weakens the assumption of metric regularity therein.
\begin{proposition}\label{lemma:lem1}
Let $\varphi:\mathbb{X}\rightarrow \mathbb{Y}$ be a mapping between two Banach spaces and $\bar x\in \mathbb{X}$. Suppose that $\varphi$ is Fr\'{e}chet differentiable at $\bar x$ and that there exist $\kappa,r>0$ such that
\begin{empheq}[box =\mybluebox]{equation}\label{3.20}
  \mathbf{d}(\bar x, \varphi^{-1}(y))\leq\kappa \|\varphi(\bar x)-y\|\ \ \forall y\in Y: \|y-\varphi(\bar x)\|<r.
\end{empheq}
 Then 
$\triangledown\varphi(\bar x)(\mathbb{X})$  is  a closed subspace of  $\mathbb{Y}$.
\end{proposition}

{\bf Proof.} Let  $y_0\in{\rm cl}(\triangledown\varphi(\bar x)(\mathbb{X}))$. Then we can find a sequence $\{y_k\}$ in $\triangledown\varphi(\bar x)(\mathbb{X})$ such that
\begin{equation}\label{3.21}
  y_k\rightarrow y_0\  {\rm and}  \ \|y_{k+1}-y_k\|<\frac{1}{2^k}\ {\rm for \ all}\ k.
\end{equation}
We claim that there exists a sequence $\{x_k\}$ in $X$ such that
\begin{equation}\label{3.22}
  \|x_{k+1}-x_k\|<\frac{3\kappa}{2^k}\ {\rm and} \ \|y_{k}-\triangledown\varphi(\bar x)(x_k)\|<\frac{1}{2^k} \ {\rm for \ all}\ k.
\end{equation}
Granting this, one has $\{x_k\}$ is a Cauchy sequence in $\mathbb{X}$ that converges to some $x_0\in \mathbb{X}$ and it follows from \eqref{3.22} that $\triangledown\varphi(\bar x)(x_k)\rightarrow y_0$, which gives $\triangledown\varphi(\bar x)(x_0)=y_0$.

We define $x_k$ iteratively. Let $x_1\in X$ be such that $\triangledown\varphi(\bar x)(x_1)=y_1$. Suppose $x_1,\cdots, x_k$ have been given to satisfy \eqref{3.22} and construct $x_{k+1}$ as follows.

Choose $u_{k+1}\in\triangledown\varphi(\bar x)^{-1}(y_{k+1})-x_k$. Let $\varepsilon>0$ be sufficiently small such that
\begin{equation}\label{3.23}
  3\varepsilon\kappa<\frac{1}{2^2}, \varepsilon<\frac{1}{2^{k+2}} \  \  {\rm and} \ \  \varepsilon\|u_{k+1}\|<\frac{1}{2^{k+2}}.
\end{equation}
Since $\varphi$ is Fr\'{e}chet differentiable at $\bar x$, then there exists $\delta\in (0,r)$ such that
\begin{equation}\label{3.24}
  \|\varphi(\bar x+w)-\varphi(\bar x)-\triangledown\varphi(\bar x)(w)\|<\varepsilon\|w\|\ \ \forall w: \|w\|<\delta.
\end{equation}
Take $t_k>0$ sufficiently small such that
$$
\max\Big\{\frac{3t_k}{2^k}, \frac{3\kappa}{2^k}t_k, t_k\|u_{k+1}\|\Big\}<\delta.
$$
Then by \eqref{3.23} and \eqref{3.24}, one has
\begin{eqnarray*}
&&\|\varphi(\bar x+t_ku_{k+1})-\varphi(\bar x)\|\\
&\leq&\|\varphi(\bar x+t_ku_{k+1})-\varphi(\bar x)-\triangledown\varphi(\bar x)(t_ku_{k+1})\|+\|t_k\triangledown\varphi(\bar x)(u_{k+1})\|\\
&<&\varepsilon t_k\|u_{k+1}\|+t_k(\|y_{k+1}-y_k\|+\|y_k-\triangledown\varphi(\bar x)(x_{k})\|)\\
&<&t_k\big(\frac{1}{2^{k+2}}+\frac{1}{2^{k}}+\frac{1}{2^{k}}\big)\\
&<&\frac{3t_k}{2^{k}}<\delta<r.
\end{eqnarray*}
This and \eqref{3.20} imply that
$$
\mathbf{d}(\bar x, \varphi^{-1}(\varphi(\bar x+t_ku_{k+1})))\leq\kappa\|\varphi(\bar x+t_ku_{k+1})-\varphi(\bar x)\|<\frac{3\kappa}{2^{k}}t_k
$$
and thus there is $w_k\in\varphi^{-1}(\varphi(\bar x+t_ku_{k+1}))$ such that
$$
\|\bar x-w_k\|<\frac{3\kappa}{2^{k}}t_k.
$$
Let $v_k:=\frac{w_k-\bar x}{t_k}$ and $x_{k+1}:=x_k+v_k$. Then
$$
\|x_{k+1}-x_k\|=\|v_k\|<\frac{3\kappa}{2^{k}}.
$$
To complete the proof, it remains to show that
\begin{equation}\label{3.25}
 \|y_{k+1}-\triangledown\varphi(\bar x)(x_{k+1})\|<\frac{1}{2^k}.
\end{equation}

Combining \eqref{3.23} with \eqref{3.24}, one has
\begin{equation}\label{3.26}
   \|\varphi(\bar x+t_kv_k)-\varphi(\bar x)-\triangledown\varphi(\bar x)(t_kv_k)\|<\varepsilon\|t_kv_k\|<\varepsilon t_k\frac{3\kappa}{2^{k}}<\frac{t_k}{2^{k+2}}
\end{equation}
and
\begin{equation}\label{3.27}
   \|\varphi(\bar x+t_ku_{k+1})-\varphi(\bar x)-\triangledown\varphi(\bar x)(t_ku_{k+1})\|<\varepsilon\|t_ku_{k+1}\|<\frac{t_k}{2^{k+2}}.
\end{equation}
Note that $\varphi(\bar x+t_kv_k)=\varphi(\bar x+t_ku_{k+1})$ and it follows from \eqref{3.26} and \eqref{3.27} that
$$
\|\triangledown\varphi(\bar x)(v_k)-\triangledown\varphi(\bar x)(u_{k+1})\|<\frac{1}{2^{k+2}}+\frac{1}{2^{k+2}}=\frac{1}{2^{k+1}}.
$$
This means that
$$
\|y_{k+1}-\triangledown\varphi(\bar x)(x_{k+1})\|=\|\triangledown\varphi(\bar x)(u_{k+1}+x_k)-\triangledown\varphi(\bar x)(x_{k}+v_k)\|<\frac{1}{2^{k+1}}.
$$
The proof is complete.\hfill$\Box$\\[1pt]

As an appropriate substitute of convexity, we consider a tangential concept introduced in 
\cite[Definition 2.1]{Sh1} under the name of $o(p)$-convexity and known in the literature  as  the Shapiro $p$-order   contact property. 
\vskip 2mm
\noindent{\bf Definition 3.2.} {\it Let $A$ be a closed subset of $\mathbb{X}$ and $p\in \N$. Recall that $A$ is said to have the $p$-order Shapiro  contact  property at $a\in A$, if for any $\varepsilon>0$ there exists $\delta>0$ such that}
    \begin{equation}\label{2-02}
      \mathbf{d}(x-u, \mathbf{T}^\mathbf{B}(A,u))\leq\varepsilon\|x-u\|^p\ \ \forall x,u\in A\cap \mathbf{B}(a,\delta).
    \end{equation}

{\it In the remainder of this paper, we will use the terminology Shapiro  first order contact property in place of $1$-order Shapiro  contact} property.\\[1pt]

In  \cite{ADT}, Aussel, Daniilids and Thibault established the links between the Shapiro  first order contact property and the notions of subsmoothness and semismoothness. We refer the reader to \cite[Theorem 3.16]{ADT} and to Thibault's book \cite[subsection 8.3.2]{thibault} for more details. In 2019 Shen, Yao and Zheng \cite{SYZ} studied the Shapiro property and  the C-Shapiro property, and showed that the Shapiro property is an extension of convexity and smoothness (see \cite[Propositions 3.8 and 3.9]{SYZ}). Further, they considered the Shapiro property of a general multifunction and used it to study calmness for closed multifunctions. Recently, inspired by the Shapiro property of a multifunction, the authors \cite{WZ2022} consider the epigraphical Shapiro property of a function; that is,\\[1pt]

\noindent{\bf Definition 3.3.} {\it Let $\varphi:\mathbb{X}\rightarrow \mathbb{R}\cup
\{+\infty\}$ be a proper function and $\bar x\in{\rm dom}(\varphi)$. We say that $\varphi$ has the epigraphical  Shapiro first order contact property  at $\bar x$, if the epigraph ${\rm epi}(\varphi)$ has the Shapiro first order contact property at $(\bar x,\varphi(\bar x))$.}\\[1pt]

The following proposition is on the epigraphical Shapiro property which was proved in \cite[Propositon 2.1]{WZ2022}.

\begin{proposition}\label{proposition:pro2}
Let $\varphi:\mathbb{X}\rightarrow \mathbb{R}\cup \{+\infty\}$ be a proper extended-real-valued lower semicontinuous function and $\bar x\in{\rm dom}(\varphi)$. Consider the following statements:
\begin{itemize}
\item[\rm(i)] $\varphi$ has   the epigraphical Shapiro first order contact property  at $\bar x$;
\item[\rm(ii)] for any $\varepsilon>0$ there exists $\delta>0$ such that
    \begin{equation}\label{2-03}
      \mathbf{d}((x-u,\varphi(x)-\varphi(u)), \mathbf{T}^\mathbf{B}({\rm epi}(\varphi),(u,\varphi(u))))\leq\varepsilon(\|x-u\|+|\varphi(x)-\varphi(u)|)
    \end{equation}
    holds for all $x,u\in \mathbf{B}_{\varphi}(\bar x,\delta)$;
\item[\rm(iii)] for any $\varepsilon>0$ there exists $\delta>0$ such that
    \begin{equation}\label{2-04}
   \varphi_H'(u;x-u)\leq\varphi(x)-\varphi(u)+\varepsilon(\|x-u\|+|\varphi(x)-\varphi(u)|)
    \end{equation}
    holds for all $x,u\in \mathbf{B}_{\varphi}(\bar x,\delta)$,
\end{itemize}
where $\mathbf{B}_{\varphi}(\bar x,\delta):=\{x\in \mathbf{B}(\bar x,\delta): |\varphi(x)-\varphi(\bar x)|<\delta\}$.

Then {\rm(i)}$\Rightarrow${\rm(ii)}$\Leftarrow${\rm(iii)}. Further, assume that $\varphi$ is continuous around $\bar x$. Then {\rm(i)}$\Leftrightarrow${\rm(ii)}.
\end{proposition}

\noindent{\it Remark 3.1.} It is noted that the implication of (ii)$\Rightarrow$(iii) or (ii)$\Rightarrow$(i) in \cref{proposition:pro2} 
 may not be necessarily true. The readers are invited to consult the counterexample in \cite[Remark 2.1]{WZ2022} for more details.\\[1pt]

The following result was proved by Shapiro and Al-Khayyal in \cite{Sh2}.
\begin{proposition}\label{proposition:pro3}
Let $\varphi:\mathbb{X}\rightarrow \mathbb{Y}$ be a continuously differentiable mapping and $A$ be a closed convex cone of $\mathbb{Y}$. Suppose that $\bar x\in \varphi^{-1}(A)$ satisfies the following Robinson qualification:
\begin{empheq}[box =\mybluebox]{equation*}
0\in{\rm int}(\varphi(\bar x)+\triangledown\varphi(\bar x)X-A).
\end{empheq}
Then $\varphi^{-1}(A)$ has the Shapiro first order contact property at $\bar x$.
\end{proposition}

The following proposition improves \cref{proposition:pro3}  
 by weakening the  Robinson qualification. This proposition is a key tool to prove main results in the paper.
\begin{proposition}\label{proposition:pro4}
Let $\varphi:\mathbb{X}\rightarrow \mathbb{Y}$ be a continuously differentiable mapping, $A$ be a closed subset of $\mathbb{Y}$ and $\bar x\in \varphi^{-1}(A)$. Suppose that $A$ has the Shapiro first order contact property at $\varphi(\bar x)$  and $\varphi$ is metrically regular around $\bar x$.
Then $\varphi^{-1}(A)$ has the Shapiro first order contact property at $\bar x$.
\end{proposition}

To prove  \cref{proposition:pro4}, we need the following two lemmata which are of independent interest.

\begin{lemma}\label{lemma:lem2}
Let $\varphi:\mathbb{X}\rightarrow \mathbb{Y}$ and $\bar x\in \mathbb{X}$. Suppose that $\varphi$ is Fr\'echet differentiable at $\bar x$ and  image $\triangledown\varphi(\bar x)(\mathbb{X})$ is of the second category. Then there exists $\mu_0>0$ such that
\begin{empheq}[box =\mybluebox]{equation}\label{3-13}
  \mathbf{d}(u,\triangledown\varphi(\bar x)^{-1}(v))\leq\mu_0\|\triangledown\varphi(\bar x)(u)-v\|\ \ \forall (u,v)\in  \mathbb{X} \times \mathbb{Y}.
\end{empheq}
Further, assume that $\varphi$ is continuously differentiable at $\bar x$. Then there exist $\mu,\delta>0$ such that
\begin{empheq}[box =\mybluebox]{equation}\label{4.13}
  \mathbf{d}(u,\triangledown\varphi(x)^{-1}(v))\leq\mu\|\triangledown\varphi(x)(u)-v\|\ \ \forall (x,u,v)\in \mathbf{B}(\bar x,\delta)\times \mathbb{X}\times \mathbb{Y}.
\end{empheq}
\end{lemma}

{\bf Proof.}  Since $\triangledown\varphi(\bar x)(\mathbb{X})$ is of the second category, it follows from the  open  mapping theorem that there exists $l>0$ such that
\begin{equation}\label{3-17}
  2l\mathbf{B}_{Y}\subseteq\triangledown\varphi(\bar x)(\mathbf{B}_\mathbb{X}).
\end{equation}
Let $(u,v)\in \mathbb{X}\times \mathbb{Y}$. By virtue of \eqref{3-17}, one has
\begin{equation}\label{4.18a}
  \triangledown\varphi(\bar x)(u)-v\in \|\triangledown\varphi(\bar x)(u)-v\|\mathbf{B}_\mathbb{Y}\subseteq\triangledown\varphi(\bar x)\Big(\frac{\|\triangledown\varphi(\bar x)(u)-v\|}{2l}\mathbf{B}_\mathbb{X}\Big).
\end{equation}
Note that $\triangledown\varphi(\bar x)^{-1}(\triangledown\varphi(\bar x)(u)-v)=u-\triangledown\varphi(\bar x)^{-1}(v)$ by the linearity of $\triangledown\varphi(\bar x)$ and then \eqref{4.18a} gives that
\begin{eqnarray*}
\mathbf{d}(u, \triangledown\varphi(\bar x)^{-1}(v))&=&\mathbf{d}(0, u-\triangledown\varphi(\bar x)^{-1}(v))\\
&=&\mathbf{d}(0, \triangledown\varphi(\bar x)^{-1}(\triangledown\varphi(\bar x)(u)-v))\\
&\leq&\frac{1}{2l}\|\triangledown\varphi(\bar x)(u)-v\|.
\end{eqnarray*}
This means that \eqref{3-13} holds with $\mu_0:=\frac{1}{2l}$.

To prove \eqref{4.13}, we next show that there exists $\delta>0$ such that
\begin{equation}\label{4.13a}
  l\mathbf{B}_{Y}\subseteq\triangledown\varphi(x)(\mathbf{B}_\mathbb{X})\ \ \forall x\in \mathbf{B}(\bar x,\delta).
\end{equation}

Noting that $\varphi$ is continuously differentiable at $\bar x$, then there exists $\delta>0$ such that $\varphi$ is Fr\'{e}chet differentiable on $\mathbf{B}(\bar x,\delta)$ and
\begin{equation}\label{4.15}
  \|\triangledown\varphi(x)-\triangledown\varphi(\bar x)\|< \frac{l}{2}\ \ \forall x\in \mathbf{B}(\bar x,\delta)
\end{equation}
(taking a smaller $\delta$ if necessary). Then by \eqref{3-17} and \eqref{4.15}, one has
$$
2l\mathbf{B}_\mathbb{Y}\subseteq(\triangledown\varphi(x)+\triangledown\varphi(\bar x)-\triangledown\varphi(x))(\mathbf{B}_\mathbb{X})\subseteq \triangledown\varphi(x)(\mathbf{B}_\mathbb{X})+\frac{l}{2}\mathbf{B}_\mathbb{Y}\ \ \forall x\in \mathbf{B}(\bar x,\delta).
$$
By the R\"{a}dstrom cancellation lemma \cite[Lemma 2.3]{Ra}, this gives that
\begin{equation}\label{4.16}
  \frac{3l}{2}\mathbf{B}_\mathbb{Y}\subseteq{\rm cl}(\triangledown\varphi(x)(\mathbf{B}_\mathbb{X}))\ \ \forall x\in \mathbf{B}(\bar x,\delta).
\end{equation}
Applying \cite[P.183, Theorem A.1]{Jameson} again gives that $\triangledown\varphi(x)(\mathbf{B}_\mathbb{X})$ and ${\rm cl}(\triangledown\varphi(x)(\mathbf{B}_\mathbb{X}))$ have the same interior and thus that \eqref{4.13a} holds by \eqref{4.16}. The proof is complete.\hfill$\Box$


\begin{lemma}\label{lemma:lem0}
Let $\varphi:\mathbb{X}\rightarrow \mathbb{Y}$ be a mapping and $A$ be a closed subset of $\mathbb{Y}$.  Suppose that $\bar x\in \varphi^{-1}(A)$ is such that $\varphi$ is Fr\'echet differentiable at $\bar x$ and metrically regular around $\bar x$. Then there exists $\delta>0$ such that
\begin{empheq}[box =\mybluebox]{equation}\label{4.a}
 \mathbf{T}^\mathbf{B}(\varphi^{-1}(A), x)=\triangledown\varphi(x)^{-1}(\mathbf{T}^\mathbf{B}(A,\varphi(x)))\ \ \forall x\in \mathbf{B}(\bar x,\delta)\cap\varphi^{-1}(A).
\end{empheq}
\end{lemma}

{\bf Proof.} Since $\varphi$ is metrically  regular around $\bar x$, there exist $\kappa>0$ along with neighborhoods $U$ of $\bar x$ and $V$ of $\varphi(\bar x)$ such that
\begin{equation}\label{4.17}
  \mathbf{d}(x,\varphi^{-1}(y))\leq\kappa\|\varphi(x)-y\|\ \ \forall (x,y)\in U\times V.
\end{equation}
Take $\delta>0$ sufficiently small such that
\begin{equation}\label{4.18}
\mathbf{B}(\bar x,\delta)\subseteq U\ \ {\rm and} \ \  \varphi(\mathbf{B}(\bar x,\delta))\subseteq V.
\end{equation}
Let $x\in \mathbf{B}(\bar x,\delta)\cap\varphi^{-1}(A)$. Choose any $u\in\triangledown\varphi(x)^{-1}(\mathbf{T}^\mathbf{B}(A,\varphi(x)))$. Then $\triangledown\varphi(x)(u)\in \mathbf{T}^\mathbf{B}(A,\varphi(x))$ and thus there exist $t_k\rightarrow 0^+$ and $v_k\rightarrow \triangledown\varphi(x)(u)$ such that
$$
\varphi(x)+t_kv_k\in A\ \ {\rm for\ all} \ k.
$$
By virtue of \eqref{4.17} and \eqref{4.18}, for any $k$ sufficiently large, one has
$$
\mathbf{d}(x+t_ku, \varphi^{-1}(\varphi(x)+t_kv_k))\leq\kappa\|\varphi(x+t_ku)-\varphi(x)-t_kv_k\|
$$
and consequently there is $x_k\in\varphi^{-1}(\varphi(x)+t_kv_k)\subseteq\varphi^{-1}(A)$ such that
\begin{equation}\label{4.19}
  \|x+t_ku-x_k\|<2\kappa\|\varphi(x+t_ku)-\varphi(x)-t_kv_k\|.
\end{equation}
Denote $u_k:=\frac{x_k-x}{t_k}$ for all $k$. Note that
$$
\varphi(x+t_ku)=\varphi(x)+\triangledown\varphi(x)(t_ku)+o(t_k)
$$
and then \eqref{4.19} gives that
$$
\|u_k-u\|<2\kappa\left\|\triangledown\varphi(x)(u)-v_k+\frac{o(t_k)}{t_k}\right\|,
$$
which implies that $u_k\rightarrow u$ as $v_k\rightarrow \triangledown\varphi(x)(u)$. Noting that $x+t_ku_k=x_k\in\varphi^{-1}(A)$, it follows that $u\in \mathbf{T}^\mathbf{B}(\varphi^{-1}(A), x)$ and thus
\begin{equation}\label{4.20}
 \mathbf{T}^\mathbf{B}(\varphi^{-1}(A), x)\supseteq\triangledown\varphi(x)^{-1}(\mathbf{T}^\mathbf{B}(A,\varphi(x))).
\end{equation}

On the other hand, let $h\in \mathbf{T}^\mathbf{B}(\varphi^{-1}(A), x)$. Then there exist $t_k\rightarrow 0^+$ and $h_k\rightarrow h$ such that $x+t_kh_k\in \varphi^{-1}(A)$ for all $k$. Using the Fr\'{e}chet differentiability of $\varphi$ at $x$, one has
$$
\varphi(x+t_kh_k)=\varphi(x)+\triangledown\varphi(x)(t_kh_k)+o(\|t_kh_k\|).
$$
Denote
$$
v_k:=\triangledown\varphi(x)(h_k)+\frac{o(\|t_kh_k\|)}{t_k}\ \ {\rm for\ all} \ k.
$$
Then $v_k\rightarrow\triangledown\varphi(x)(h)$ and $\varphi(x)+t_kv_k=\varphi(x+t_kh_k)\in A$, which implies that $\triangledown\varphi(x)(h)\in \mathbf{T}^\mathbf{B}(A,\varphi(x))$. Hence
$$
\mathbf{T}^\mathbf{B}(\varphi^{-1}(A), x)\subseteq\triangledown\varphi(x)^{-1}(\mathbf{T}^\mathbf{B}(A,\varphi(x)))
$$
and \eqref{4.a} holds by \eqref{4.20}. The proof is complete.\hfill$\Box$\\[1pt]

{\it Proof of  \cref{proposition:pro4}}. 
 Since $\varphi$ is continuously differentiable and metrically regular around $\bar x$, it follows from \cite[Theorem 1.57]{M2} that $\triangledown\varphi(\bar x)$ is surjective. By virtue of  \cref{lemma:lem2} and \cref{lemma:lem0} ,
there exist $\mu,\delta>0$ such that \eqref{4.13} and \eqref{4.a} hold.

Let $\varepsilon>0$.  Since $\varphi$ is continuously differentiable at $\bar x$, there exists $\delta_1\in(0,\delta)$ such that
\begin{equation}\label{4.24}
  \|\varphi(x)-\varphi(u)-\triangledown\varphi(u)(x-u)\|<\frac{\varepsilon}{2\mu}\|x-u\|\ \ \forall x,u\in \mathbf{B}(\bar x,\delta_1)
\end{equation}
and
\begin{equation}\label{4.25}
  \|\varphi(x)-\varphi(u)\|\leq K\|x-u\|\ \ \forall x,u\in \mathbf{B}(\bar x,\delta_1)
\end{equation}
where $K:=\|\triangledown\varphi(\bar x)\|+1$.

Noting $A$ has the Shapiro first order contact property at $\varphi(\bar x)$, it follows that there exists $\delta_2\in (0,\delta_1)$ such that
\begin{equation}\label{4.26}
  \mathbf{d}(y-v, \mathbf{T}^\mathbf{B}(A,v))\leq\frac{\varepsilon}{2K\mu}\|y-v\|\ \ \forall y,v\in A\cap \mathbf{B}(\varphi(\bar x),\delta_2).
\end{equation}

Take $\delta_3\in(0,\delta_2)$ such that $\varphi(\mathbf{B}(\bar x,\delta_3))\subseteq \mathbf{B}(\varphi(\bar x),\delta_2)$. Let $x,u\in \varphi^{-1}(A)\cap \mathbf{B}(\bar x,\delta_3)$. Take a sequence $\{v_n\}$ in $\mathbf{T}^\mathbf{B}(A,\varphi(u))$ such that
$$
\|\varphi(x)-\varphi(u)-v_n\|\rightarrow \mathbf{d}(\varphi(x)-\varphi(u), \mathbf{T}^\mathbf{B}(A,\varphi(u))).
$$
By virtue of \eqref{4.13} and \eqref{4.a}, one has
\begin{eqnarray*}
\mathbf{d}(x-u, \mathbf{T}^\mathbf{B}(\varphi^{-1}(A), u))&=&\mathbf{d}(x-u,\triangledown\varphi(u)^{-1}(\mathbf{T}^\mathbf{B}(A,\varphi(u))))\\
&\leq&\mathbf{d}(x-u,\triangledown\varphi(u)^{-1}(v_n))\\
&\leq&\mu\|\triangledown\varphi(u)(x-u)-v_n\|\\
&\leq&\mu\big(\|\triangledown\varphi(u)(x-u)+\varphi(u)-\varphi(x)\|+\|\varphi(x)-\varphi(u)-v_n\|\big).
\end{eqnarray*}
By letting $n\rightarrow\infty$, one gets
\begin{eqnarray*}
&&\mathbf{d}(x-u, \mathbf{T}^\mathbf{B}(\varphi^{-1}(A), u))\\
&\leq&\mu\big(\|\varphi(x)-\varphi(u)-\triangledown\varphi(u)(x-u)\|+\mathbf{d}(\varphi(x)-\varphi(u), \mathbf{T}^\mathbf{B}(A,\varphi(u)))\|\big).
\end{eqnarray*}
This together with \eqref{4.24}, \eqref{4.25} and \eqref{4.26} gives
\begin{eqnarray*}
&&\mathbf{d}(x-u, \mathbf{T}^\mathbf{B}(\varphi^{-1}(A), u))\\
&\leq&\mu\big(\|\varphi(x)-\varphi(u)-\triangledown\varphi(u)(x-u)\|+\mathbf{d}(\varphi(x)-\varphi(u), \mathbf{T}^\mathbf{B}(A,\varphi(u)))\|\big)\\
&\leq& \mu\big(\frac{\varepsilon}{2\mu}\|x-u\|+\frac{\varepsilon}{2K\mu}\|\varphi(x)-\varphi(u)\|\big)\\
&\leq&\mu\big(\frac{\varepsilon}{2\mu}\|x-u\|+\frac{\varepsilon}{2K\mu}K\|x-u\|\big)\\
&=&\varepsilon\|x-u\|.
\end{eqnarray*}
Hence $\varphi^{-1}(A)$ has the Shapiro first order contact property  at $\bar x$. The proof is complete.\hfill$\Box$\\[1pt]
\setcounter{equation}{0}
\section{Main results}
In this section, we study error bounds of the inequality defined by a general function and aim to provide primal criteria of error bounds via Bouligand tangent cones,  lower  Hadamard directional derivatives and the Hausdorff-Pompeiu excess  of subsets. Then we apply these primal results to error bounds of the  composite-convex function, and establish primal characterizations of error bounds and an accurate estimate of the error bound modulus in terms of Bouligand tangent cones, directional derivatives of convex functions and the Hausdorff-Pompeiu excess.  We first consider error bounds of the inequality defined by a general lower semicontinuous function.

Given a proper lower semicontinuous extended-real-valued function $\varphi:\mathbb{X}\rightarrow\mathbb{R}\cup\{+\infty\}$, we consider the following inequality:
\begin{empheq}[box =\mybluebox]{equation}\label{3.1}
 \varphi(x)\leq 0.
 \end{empheq}
We denote by $\mathbf{S}_{\varphi}:=\{x\in \mathbb{X}: \varphi(x)\leq 0\}$ the solution set. Recall that inequality \eqref{3.1} is said to have a local error bound at $\bar x\in \mathbf{S}_{\varphi}$, if there exist $\tau,\delta\in(0, +\infty)$ such that
\begin{empheq}[box =\mybluebox]{equation}\label{3.2}
  \mathbf{d}(x, \mathbf{S}_{\varphi})\leq \tau [\varphi(x)]_+,\ \ \forall x\in \mathbf{B}(\bar x,\delta),
\end{empheq}
where $[\varphi(x)]_+:=\max\{\varphi(x), 0\}$. We denote by
\begin{equation}\label{3.3a}
\tau(\mathbf{S}_{\varphi},\bar x):=\inf\{\tau>0: \ {\rm there\ exists} \ \delta>0\ {\rm such\ that} \ \eqref{3.2}\ {\rm holds}\}
\end{equation}
the local error bound modulus of $\mathbf{S}_{\varphi}$ at $\bar x$.

The theory of error bounds has been a subject of intense study during many years  due to its numerous applications in optimization and variational analysis. They include  for instance areas like sensitivity analysis of linear programming, convergence analysis of descent methods, the feasibility problems  and the domain of image reconstruction. This notion has been proved to have close connections with several notions in convex analysis and approximation like the  basic constraint qualification (BCQ),  the Abadie constraint qualification (ACQ), the strong  conical hull intersection property (CHIP), the linear regularity and normal property, the metric subregularity as well as the calmness of multifunctions.

The following theorem gives necessary and/or sufficient criteria of local error bounds for the inequality \eqref{3.1} in terms of the   Bouligand tangent cone, the  lower  Hadamard directional derivative and the Hausdorff-Pompeiu excess  of a set beyond another set.

\begin{theorem}\label{theorem:th1}
Let $\varphi:X\rightarrow\mathbb{R}\cup\{+\infty\}$ be a proper lower semicontinuous function and $\bar x\in \mathbf{S}_{\varphi}$.
\begin{itemize}
\item[\rm(i)] Suppose that $\mathbf{S}_{\varphi}$ has the Shapiro first order contact property around $\bar x$ and inequality \eqref{3.1} has a local error bound at $\bar x$. Then there exist $\tau, r\in(0, +\infty)$ such that
\begin{equation}\label{3.3}
  \mathbf{e}\big(\varphi'_\mathbf{H}(x,\cdot)^{-1}(-\infty, 1], \mathbf{T}^\mathbf{B}(\mathbf{S}_{\varphi}, x)\big)\leq \tau
\end{equation}
holds for all $x\in \mathbf{S}_{\varphi}\cap \mathbf{B}(\bar x, r)$ with $\varphi(x)=0$.
\item[\rm(ii)] Suppose that ${\rm bd}(\mathbf{S}_{\varphi})\subseteq\varphi^{-1}(0)$, $\varphi$ has the epigraphical Shapiro first order contact property  at $\bar x$ and that there exist $\tau, r\in(0, +\infty)$ such that \eqref{3.3} holds for all $x\in {\rm bd}(\mathbf{S}_{\varphi})\cap \mathbf{B}(\bar x, r)$. Then inequality \eqref{3.1} has a local error bound at $\bar x$.
\end{itemize}
\end{theorem}

{\bf Proof.} (i)  By the local error bound of inequality \eqref{3.1} at $\bar x$, there exist $\tau,\delta\in(0, +\infty)$ such that \eqref{3.2} holds. Note that $\mathbf{S}_{\varphi}$ has the Shapiro first order contact property around $\bar x$ and thus there exists $r\in (0,\delta)$ such that $\mathbf{S}_{\varphi}$ has the Shapiro first order contact property  on $\mathbf{S}_{\varphi}\cap \mathbf{B}(\bar x, r)$. Let $x\in \mathbf{S}_{\varphi}\cap \mathbf{B}(\bar x, r)$ with $\varphi(x)=0$. Take any $h\in X$ such that $\varphi'_\mathbf{H}(x,h)\leq 1$. We need to show that
\begin{equation}\label{3.4}
  \mathbf{d}(h,\mathbf{T}^\mathbf{B}(\mathbf{S}_{\varphi}, x))\leq \tau.
\end{equation}

Note that $\varphi'_\mathbf{H}(x,h)\leq 1$ and then there exist $t_n\rightarrow 0^+$ and $h_n\rightarrow h$ such that
\begin{equation}\label{3.5}
  \frac{\varphi(x+t_nh_n)-\varphi(x)}{t_n}<1+\frac{1}{n}.
\end{equation}
If there exists a subsequence $\{n_k\}\subseteq\mathbb{N}$ such that $\varphi(x+t_{n_k}h_{n_k})\leq0$ for all $k$, then one has $h\in \mathbf{T}^\mathbf{B}(\mathbf{S}_{\varphi},x)$ and thus \eqref{3.4} holds by $\mathbf{d}(h, \mathbf{T}^\mathbf{B}(\mathbf{S}_{\varphi},x))=0$.

We next consider the case that $\varphi(x+t_nh_n)>0$ for all $n$. Since $\mathbf{S}_{\varphi}$ has the Shapiro first order contact property  at $x$, by Proposition 2.1, for any $k\in\mathbb{N}$, there exist $r_k\in(0, +\infty)$ such that $\mathbf{B}(x, r_k)\subseteq \mathbf{B}(\bar x,\delta)$ and
\begin{equation}\label{3.6}
  \mathbf{d}(u-x, \mathbf{T}^\mathbf{B}(\mathbf{S}_{\varphi}, x))\leq \mathbf{d}(u, \mathbf{S}_{\varphi})+\frac{1}{k}\|u-x\|\ \ \forall u\in \mathbf{B}(x, r_k).
\end{equation}
Take a subsequence $\{n_k\}\subseteq\mathbb{N}$ such that $\|t_{n_k}h_{n_k}\|<\frac{1}{k}$. Then by virtue of \eqref{3.2}, \eqref{3.5} and \eqref{3.6}, one has
\begin{eqnarray*}
\mathbf{d}(x+t_{n_k}h_{n_k}-x, \mathbf{T}^\mathbf{B}(\mathbf{S}_{\varphi},x))&\leq& \mathbf{d}(x+t_{n_k}h_{n_k}, \mathbf{S}_{\varphi})+\frac{1}{k} t_{n_k}\|h_{n_k}\|\\
&\leq& \tau\varphi(x+t_{n_k}h_{n_k})+\frac{1}{k} t_{n_k}\|h_{n_k}\|\\
&<& \tau\big(1+\frac{1}{n_k}\big) t_{n_k}+\frac{1}{k} t_{n_k}\|h_{n_k}\|
\end{eqnarray*}
thanks to $\varphi(x)=0$. This implies that
\begin{eqnarray*}
\mathbf{d}(h_{n_k}, \mathbf{T}^\mathbf{B}(\mathbf{S}_{\varphi},x))\leq \tau\big(1+\frac{1}{n_k}\big)+\frac{1}{k}\|h_{n_k}\|
\end{eqnarray*}
and consequently \eqref{3.4} holds by letting $k\rightarrow\infty$.

(ii) Let $\varepsilon>0$ such that $(1+\tau)\varepsilon<1$. Since $\varphi$ has the Shapiro first order contact property at $\bar x$, then there exists $\delta_1\in(0, r)$ such that
\begin{equation}\label{3.7}
  \mathbf{d}((x-u,\alpha-\beta), \mathbf{T}^\mathbf{B}({\rm epi}(\varphi), (u,\beta)))<\varepsilon(\|x-u\|+\|\alpha-\beta\|)
\end{equation}
holds for all $(x,\alpha),(u,\beta)\in{\rm epi}(\varphi)\cap (\mathbf{B}(\bar x,\delta_1)\times (\varphi(\bar x)-\delta_1,\varphi(\bar x)+\delta_1))$.

Take $\delta\in (0, \frac{\delta_1}{2})$ and let $x\in \mathbf{B}(\bar x, \delta)\backslash \mathbf{S}_{\varphi}$. Then $\mathbf{d}(x,\mathbf{S}_{\varphi})\leq\|x-\bar x\|<\delta$. Choose any $\gamma\in (0,1)$ such that
\begin{equation}\label{3.8}
  \gamma>\max\Big\{\frac{\mathbf{d}(x,\mathbf{S}_{\varphi})}{\delta}, (1+\tau)\varepsilon\Big\}.
\end{equation}
By using \cref{lemma:lema},
there exists $a\in {\rm bd}(\mathbf{S}_{\varphi})$ such that
\begin{equation}\label{3.9}
  \gamma\|x-a\|\leq\min\{\mathbf{d}(x-a, \mathbf{T}^\mathbf{B}(\mathbf{S}_{\varphi},a)), \mathbf{d}(x,\mathbf{S}_{\varphi})\}.
\end{equation}
Then
$$
\|\bar x-a\|\leq\|\bar x-x\|+\|x-a\|<\delta+\frac{\mathbf{d}(x,\mathbf{S}_{\varphi})}{\gamma}<2\delta<\delta_1
$$
and it follows from \eqref{3.7} that
$$
\mathbf{d}((x-a,\varphi(x)-\varphi(a)), \mathbf{T}^\mathbf{B}({\rm epi}(\varphi), (a,\varphi(a))))<\varepsilon(\|x-a\|+|\varphi(x)-\varphi(a)|).
$$
Thus, there is $(h,\alpha)\in \mathbf{T}^\mathbf{B}({\rm epi}(\varphi), (a,\varphi(a)))={\rm epi}(\varphi_H'(a,\cdot))$ such that
\begin{equation}\label{3.10}
 \varphi_H'(a, h)\leq\alpha\ \ {\rm and} \ \ \|x-a-h\|+|\varphi(x)-\varphi(a)-\alpha|<\varepsilon(\|x-a\|+|\varphi(x)-\varphi(a)|).
\end{equation}

Case 1: $\alpha=0$. Then for any $\lambda>0$, one has $\varphi_H'(a, \lambda h)\leq 0 $ and then \eqref{3.3} implies that
$$
\mathbf{d}(\lambda h, \mathbf{T}^\mathbf{B}(\mathbf{S}_{\varphi},a))\leq 1.
$$
This means $h\in \mathbf{T}^\mathbf{B}(\mathbf{S}_{\varphi},a)$ by letting $\lambda\rightarrow+\infty$. By virtue of \eqref{3.9} and \eqref{3.10}, one has
\begin{eqnarray*}
\gamma\|x-a\|\leq \mathbf{d}(x-a, \mathbf{T}^\mathbf{B}(\mathbf{S}_{\varphi},a))\leq \|x-a-h\|&<&\varepsilon(\|x-a\|+|\varphi(x)-\varphi(a)|)\\
&=&\varepsilon\|x-a\|+\varepsilon\varphi(x)
\end{eqnarray*}
the equality follows by $\varphi(x)>0$ and $\varphi(a)=0$. This implies that
$$
(\gamma-\varepsilon)\mathbf{d}(x, \mathbf{S}_{\varphi})\leq(\gamma-\varepsilon)\|x-a\|<\varepsilon\varphi(x).
$$
By letting $\gamma\rightarrow 1^-$, one has
\begin{equation}\label{3.11}
\mathbf{d}(x, \mathbf{S}_{\varphi})\leq\frac{\varepsilon}{1-\varepsilon}\varphi(x).
\end{equation}

Case 2: $\alpha\not=0$. Note that $\varphi_H'(a,\cdot)$ is positively homogeneous and then by \eqref{3.10}, one has
$$
 \varphi_H'\big(a, \frac{h}{|\alpha|}\big)\leq 1.
$$
This and \eqref{3.3} imply that
$$
d\big(\frac{h}{|\alpha|}, \mathbf{T}^\mathbf{B}(\mathbf{S}_{\varphi},a)\big)\leq \tau
$$
and consequently
$$
\mathbf{d}(h, \mathbf{T}^\mathbf{B}(\mathbf{S}_{\varphi},a))\leq \tau|\alpha|.
$$
By virtue of \eqref{3.9} and \eqref{3.10}, one has
\begin{eqnarray*}
\gamma\|x-a\|&\leq& \mathbf{d}(x-a,\mathbf{T}^\mathbf{B}(\mathbf{S}_{\varphi},a))\\
&\leq& \|x-a-h\|+ \mathbf{d}(h,\mathbf{T}^\mathbf{B}(\mathbf{S}_{\varphi},a))\\
&\leq& \|x-a-h\|+\tau(|\varphi(x)-\varphi(a)-\alpha|+|\varphi(x)-\varphi(a)|)\\
&<&\varepsilon(\|x-a\|+|\varphi(x)-\varphi(a)|)+\tau(\varepsilon\|x-a\|+\varepsilon|\varphi(x)-\varphi(a)|+|\varphi(x)-\varphi(a)|)
\end{eqnarray*}
and thus
\begin{eqnarray*}
(\gamma-\varepsilon-\tau\varepsilon)\mathbf{d}(x,\mathbf{S}_{\varphi})\leq(\gamma-\varepsilon-\tau\varepsilon)\|x-a\|&<&
(\tau+\varepsilon+\tau\varepsilon)|\varphi(x)-\varphi(a)|\\
&=&(\tau+\varepsilon+\tau\varepsilon)\varphi(x)
\end{eqnarray*}
the equality holds by $\varphi(x)>0$ and $\varphi(a)=0$. By letting $\gamma\rightarrow 1^-$, one has
\begin{eqnarray*}
(1-\varepsilon-\tau\varepsilon)\mathbf{d}(x,\mathbf{S}_{\varphi})<(\tau+\varepsilon+\tau\varepsilon)\varphi(x)
\end{eqnarray*}
and consequently
\begin{eqnarray*}
\mathbf{d}(x,\mathbf{S}_{\varphi})\leq\frac{\tau+\varepsilon+\tau\varepsilon}{1-\varepsilon-\tau\varepsilon}\varphi(x).
\end{eqnarray*}
This and \eqref{3.11} imply that inequality \eqref{3.1} has the local error bound at $\bar x$ (with the constant $\max\{\frac{\varepsilon}{1-\varepsilon}, \frac{\tau+\varepsilon+\tau\varepsilon}{1-\varepsilon-\tau\varepsilon}\}$). The proof is complete.\hfill$\Box$
\\[1pt]

The following theorem, immediate from (i) and (ii) in \cref{theorem:th1}
provides an equivalent primal condition for the local error bound under the Shapiro first order contact property.

\begin{theorem}\label{theorem:th2}
Let $\bar x\in \mathbf{S}_{\varphi}$. Suppose that ${\rm bd}(\mathbf{S}_{\varphi})\subseteq\varphi^{-1}(0)$, $\varphi$ has the epigraphical Shapiro first order contact property  at $\bar x$ and that $\mathbf{S}_{\varphi}$ has the Shapiro first order contact property   around $\bar x$. Then inequality \eqref{3.1} has the local error bound at $\bar x$ if and only if
\begin{equation*}
\limsup_{x\stackrel{{\rm bd}(\mathbf{S}_{\varphi})}\longrightarrow \bar x}\mathbf{e}\big(\varphi'_\mathbf{H}(x,\cdot)^{-1}(-\infty, 1], \mathbf{T}^\mathbf{B}(\mathbf{S}_{\varphi}, x)\big)<+\infty.
\end{equation*}
Further, one has the following accurate estimate for the local error bound modulus:
\begin{equation}\label{3.12}
 \tau(\mathbf{S}_{\varphi}, \bar x)=\limsup_{x\stackrel{{\rm bd}(\mathbf{S}_{\varphi})}\longrightarrow \bar x}\mathbf{e}\big(\varphi'_\mathbf{H}(x,\cdot)^{-1}(-\infty, 1], \mathbf{T}^\mathbf{B}(\mathbf{S}_{\varphi}, x)\big).
\end{equation}
\end{theorem}

{\bf Proof.} Thanks to (i) and (ii) in  \cref{theorem:th1} ,
it suffices to prove \eqref{3.12}. We denote
$$
\alpha:=\limsup_{x\stackrel{{\rm bd}(\mathbf{S}_{\varphi})}\longrightarrow \bar x}\mathbf{e}\big(\varphi'_\mathbf{H}(x,\cdot)^{-1}(-\infty, 1], \mathbf{T}^\mathbf{B}(\mathbf{S}_{\varphi}, x)\big).
$$

We first consider the case $\tau(\mathbf{S}_{\varphi}, \bar x)<+\infty$. Let $\tau\in (0, \tau(\mathbf{S}_{\varphi}, \bar x))$. Then by using the proof of (i) in  \cref{theorem:th1},
there exists $\delta>0$ such that
$$
\mathbf{e}\big(\varphi'_\mathbf{H}(x,\cdot)^{-1}(-\infty, 1], \mathbf{T}^\mathbf{B}(\mathbf{S}_{\varphi}, x)\big)\leq \tau\ \ \forall x\in \mathbf{B}(\bar x,\delta)\cap{\rm bd}(\mathbf{S}_{\varphi}).
$$
This means that $\alpha\leq \tau$ and thus $\alpha\leq \tau(\mathbf{S}_{\varphi}, \bar x)$ by letting $\tau\downarrow \tau(\mathbf{S}_{\varphi}, \bar x)$.

Note that $\alpha<+\infty$ and for any $\nu>0$ there exists $\delta>0$ such that
$$
\mathbf{e}\big(\varphi'_\mathbf{H}(x,\cdot)^{-1}(-\infty, 1], \mathbf{T}^\mathbf{B}(\mathbf{S}_{\varphi}, x)\big)\leq \alpha+\nu \ \ \forall x\in \mathbf{B}(\bar x,\delta)\cap{\rm bd}(\mathbf{S}_{\varphi}).
$$
Applying the proof of (ii) in  \cref{theorem:th1}, 
for any $\varepsilon>0$ sufficiently small, one has
$$
\tau(\mathbf{S}_{\varphi}, \bar x)\leq \frac{(\alpha+\nu)+\varepsilon+(\alpha+\nu)\varepsilon}{1-\varepsilon-(\alpha+\nu)\varepsilon}.
$$
By letting $\varepsilon\downarrow 0$, one has
$$
\tau(\mathbf{S}_{\varphi}, \bar x)\leq \alpha+\nu\  \ \forall \nu>0,
$$
and consequently $\tau(\mathbf{S}_{\varphi}, \bar x)\leq \alpha$ by letting $\nu\downarrow 0$. Hence \eqref{3.12} holds.

We next consider the case $\tau(\mathbf{S}_{\varphi}, \bar x)=+\infty$. We claim that $\alpha=+\infty$ (otherwise, $\alpha<+\infty$ and by using the proof of (ii) in \cref{theorem:th1}
 again, one can obtain that $\tau(\mathbf{S}_{\varphi}, \bar x)<+\infty$, a contradiction). The proof is complete.\hfill$\Box$
\\[1pt]

\noindent{\it Remark 4.1.} For the case that $\varphi$ is convex, the Shapiro first order contact property  holds automatically and thus \cref{theorem:th2}
reduces to \cite[Theorem 5.3]{WZ2021}, which means that  \cref{theorem:th2} 
is an extension of \cite[Theorem 5.3]{WZ2021} from the convex case to the non-convex one.\\[1pt]

It is noted that \eqref{3.3} is a key inequality to characterize the local error bound of inequality \eqref{3.1}. We are now in a position to give a characterization of \eqref{3.3} via the following proposition.
\begin{proposition}\label{proposition:pro5}
Let $\tau>0$ and $x\in \mathbf{S}_{\varphi}$. Then \eqref{3.3} holds if and only if
\begin{empheq}[box =\mybluebox]{equation}\label{3.14}
  \mathbf{d}(h, \mathbf{T}^\mathbf{B}(\mathbf{S}_{\varphi},x))\leq \tau \max\{\varphi_H'(x,h), 0\}\ \ \forall h\in X.
\end{empheq}
\end{proposition}

{\bf Proof.} The necessity part. Suppose that \eqref{3.3} holds. Let $h\in X$ be such that $\varphi_H'(x,h)\leq 0$. Then for any $\lambda>0$, one has $\varphi_H'(x, \lambda h)\leq 0$ and it follows from \eqref{3.3} that
$$
\mathbf{d}(\lambda h, \mathbf{T}^\mathbf{B}(\mathbf{S}_{\varphi},x))\leq \tau
$$
This implies that
$$
\mathbf{d}(h, \mathbf{T}^\mathbf{B}(\mathbf{S}_{\varphi},x))\leq \frac{\tau}{\lambda} \ \ {\rm for\ any } \ \lambda>0.
$$
By letting $\lambda\rightarrow +\infty$, one has $h\in \mathbf{T}^\mathbf{B}(\mathbf{S}_{\varphi},x)$ as $\mathbf{T}^\mathbf{B}(\mathbf{S}_{\varphi},x)$ is closed, which implies that \eqref{3.14} holds.

Let $h\in X$ be such that $\varphi_H'(x,h)> 0$. Note that $\varphi_H'(a,\cdot)$ is positively homogeneous and then
$$
\varphi_H'\left(\mathbf{T}^\mathbf{B}(x, \frac{h}{\varphi_H'(x,h)}\right)=1.
$$
By virtue of \eqref{3.3}, one has
$$
 d\left(\mathbf{T}^\mathbf{B}(\frac{h}{\varphi_H'(x,h)}h, \mathbf{T}^\mathbf{B}(\mathbf{S}_{\varphi},x)\right)\leq \tau
$$
and consequently
$$
 \mathbf{d}(h, \mathbf{T}^\mathbf{B}(\mathbf{S}_{\varphi},x))\leq \tau\varphi_H'(x,h),
$$
which implies that \eqref{3.14} holds.

The sufficiency part. Suppose that \eqref{3.14} holds. Let $h\in X$ be such that $\varphi_H'(x,h)\leq 1$. Then \eqref{3.14} implies that
$$
\mathbf{d}(h,\mathbf{T}^\mathbf{B}(\mathbf{S}_{\varphi},x))\leq\tau\max\{\varphi_H'(x,h), 0\}\leq \tau.
$$
Hence \eqref{3.3} holds. The proof is complete.\hfill$\Box$
\\[1pt]

The following theorem follows immediately from  \cref{theorem:th1} 
 and  \cref{proposition:pro5}. 

\begin{theorem}\label{theorem:th3}
Let $\bar x\in \mathbf{S}_{\varphi}$ and $\tau>0$. Suppose that ${\rm bd}(\mathbf{S}_{\varphi})\subseteq\varphi^{-1}(0)$, $\varphi$ has the epigraphical Shapiro first order contact property  at $\bar x$ and that $\mathbf{S}_{\varphi}$ has the Shapiro first order contact property around $\bar x$. Then the following statements are equivalent:
\begin{itemize}
\item[\rm(i)] Inequality \eqref{3.1} has the local error bound at $\bar x$ with constant $\tau>0$;
\item[\rm(ii)] There exists $\delta>0$ such that \eqref{3.3} holds for all $x\in {\rm bd}(\mathbf{S}_{\varphi})\cap \mathbf{B}(\bar x,\delta)$;
\item[\rm(ii)] There exists $\delta>0$ such that \eqref{3.14} holds for all $x\in {\rm bd}(\mathbf{S}_{\varphi})\cap \mathbf{B}(\bar x,\delta)$.
\end{itemize}
\end{theorem}

Given $x\in X$, we can consider the inequality $\varphi_H'(x,\cdot)\leq 0$ defined by the  lower  Hadamard directional derivative $\varphi_H'(x,\cdot)$ and study global error bounds
of the inequality $\varphi_H'(x,\cdot)\leq 0$; that is, there exists $\tau>0$ such that
\begin{equation}\label{3.15}
  \mathbf{d}(h,S_{\varphi_H'(x,\cdot)})\leq \tau\max\{\varphi_H'(x,h), 0\}\ \ {\rm for\ all} \ h\in X,
\end{equation}
where $S_{\varphi_H'(x,\cdot)}:=\{u\in X:\varphi_H'(x,u)\leq 0\}$. \\[1pt]

The following corollary, immediate from  \cref{theorem:th1}, shows that the local error bound of $\mathbf{S}_{\varphi}$ at $\bar x$ is, to some degree, equivalent to the global error bound of the inequality $\varphi_H'(x,\cdot)\leq 0$ for all $x$ close to $\bar x$ with the same constant. This corollary is inspired from \cite[Theorem 4.4]{WZ2018} which studies the convex inequality defined by the Clarke directional derivative of a locally Lipschitz function and its global error bound.

\begin{corollary}
Let $\bar x\in \mathbf{S}_{\varphi}$ and $\tau>0$. Assume that ${\rm bd}(\mathbf{S}_{\varphi})\subseteq\varphi^{-1}(0)$, $\varphi$ has the epigraphical Shapiro first order contact property at $\bar x$ and that $\mathbf{S}_{\varphi}$ has the Shapiro first order contact property  around $\bar x$.
\begin{itemize}
\item[\rm(i)] Suppose that inequality \eqref{3.1} has the local error bound at $\bar x$. Then there exists $\delta>0$ such that for any $x\in {\rm bd}(\mathbf{S}_{\varphi})\cap \mathbf{B}(\bar x,\delta)$, the inequality $\varphi_H'(x,\cdot)\leq 0$ has the global error bound with the same constant.
\item[\rm(ii)] Suppose that there exists a neighborhood $U$ of $\bar x$ such that $\ker\varphi_H'(x,\cdot)\subseteq \mathbf{T}^\mathbf{B}(\mathbf{S}_{\varphi}, x)$ for all $x\in U\cap{\rm bd}(\mathbf{S}_{\varphi})$. Then inequality \eqref{3.1} has the local error bound at $\bar x$ if and only if there exists $\delta>0$ such that for any $x\in {\rm bd}(\mathbf{S}_{\varphi})\cap \mathbf{B}(\bar x,\delta)$, the inequality $\varphi_H'(x,\cdot)\leq 0$ has the global error bound with the same constant.
\end{itemize}
\end{corollary}

{\bf Proof.} By the definition of the  lower  Hadamard directional derivative, for any $x\in {\rm bd}(\mathbf{S}_{\varphi})$, one can verify that
$$
\{h\in X:\varphi_H'(x,h)<0\}\subseteq \mathbf{T}^\mathbf{B}(\mathbf{S}_{\varphi}, x)\subseteq \{h\in X:\varphi_H'(x,h)\leq 0\}
$$
and thus the conclusions follow from  \cref{theorem:th3}.
The proof is complete.\hfill$\Box$
\\[1pt]



Now, we are in a position to study error bounds of a composite-convex function. We apply main results obtained above to establish primal characterizations of error bounds and give an accurate estimate of the error bound modulus in terms of Bouligand tangent cones, directional derivatives of convex functions and the Hausdorff-Pompeiu excess.\\[1pt]

{\it Throughout the rest of this section, we always assume that that $\mathbb{Y}$ is a Banach space, $g:\mathbb{X}\rightarrow \mathbb{Y}$ is a continuously differentiable mapping and that $f:\mathbb{Y}\rightarrow\mathbb{R}\cup\{+\infty\}$ is a proper lower semicontinuous and convex function.}\\[1pt]

We consider the following composite-convex inequality:
\begin{empheq}[box =\mybluebox]{equation}\label{3.16}
  (f\circ g)(x)\leq 0.
\end{empheq}
We denote by $\mathcal{S}:=\{x\in X:  (f\circ g)(x)\leq 0\}$ the solution set of \eqref{3.16}. For any given $\bar x\in\mathcal{S}$,  we denote by $\tau(\mathcal{S},\bar x)$, defined as said in \eqref{3.3a}, the local error bound modulus of $\mathcal{S}$ at $\bar x$.

\begin{theorem}\label{theorem:th4}
Denote $S_f:=\{y\in Y: f(y)\leq 0\}$ and assume that ${\rm bd}(S_f)\subseteq f^{-1}(0)$. Let $\bar x\in\mathcal{S}$ be such that $g(\bar x)\in{\rm int}({\rm dom}f)$ and $g$ is metrically regular around $\bar x$.  Then composite-convex inequality \eqref{3.16} has a local error bound at $\bar x$ if and only if
\begin{empheq}[box =\mybluebox]{equation*}
\limsup_{x\stackrel{{\rm bd}(\mathcal{S})}\longrightarrow \bar x}\;\mathbf{e}\big(\{h\in X: \mathbf{d}^+f(g(x), \triangledown g(x)(h))\leq 1\},  \triangledown g(x)^{-1}(\mathbf{T}^\mathbf{B}(S_{f}, g(x))\big)<+\infty.
\end{empheq}
Further, one has the following accurate estimate for the local error bound modulus:
\begin{empheq}[box =\mybluebox]{equation*}
 \tau(\mathcal{S}, \bar x)=\limsup_{x\stackrel{{\rm bd}(\mathcal{S})}\longrightarrow \bar x}\mathbf{e}\big(\{h\in X: \mathbf{d}^+f(g(x), \triangledown g(x)(h))\leq 1\},  \triangledown g(x)^{-1}(\mathbf{T}^\mathbf{B}(S_{f}, g(x))\big).
\end{empheq}
\end{theorem}

To prove  \cref{theorem:th4}, 
we need the following lemma which is of some independent interest.

\begin{lemma}\label{lemma:lem10}
Let $\varphi:=f\circ g$ and $\bar x\in X$ be such that $g(\bar x)\in{\rm int}({\rm dom}f)$. Then there exists $\delta>0$ such that
 \begin{equation}\label{3.17}
\varphi_H'(x,h)=\mathbf{d}^+f(g(x),\triangledown g(x)(h))\ \ \forall h\in X
\end{equation}
holds for all $x\in \mathbf{B}(\bar x,\delta)$.
\end{lemma}

{\bf Proof.} Note that $g(\bar x)\in{\rm int}({\rm dom}f)$ and then \cite[Proposition 1.6]{Ph} implies that $f$ is locally Lipschtizian around $g(\bar x)$; that is, there exist $L,r>0$ such that
\begin{equation}\label{3.18}
|f(y_1)-f(y_2)|\leq L\|y_1-y_2\|\ \ \forall y_1,y_2\in \mathbf{B}(g(\bar x), r).
\end{equation}
Using the continuity of $g$, there is $\delta>0$ such that
\begin{equation}\label{3.19}
  g(\mathbf{B}(\bar x,\delta))\subseteq \mathbf{B}(g(\bar x), r).
\end{equation}
Let $x\in \mathbf{B}(\bar x,\delta)$ and take any $h\in X$. Then for any $t>0$ sufficiently small and $h'$ close to $h$, one has
\begin{eqnarray*}
&&\frac{1}{t}[\varphi(x,th')-\varphi(x)]=\frac{1}{t}[f(g(x,th'))-f(g(x))]\\
&=&\frac{1}{t}[f(g(x)+\triangledown g(x)(th')+o(t))-f(g(x))]\\
&=&\frac{1}{t}[f(g(x)+\triangledown g(x)(th')+o(t))-f(g(x)+\triangledown g(x)(th))]\\
&&+\frac{1}{t}[f(g(x)+\triangledown g(x)(th))-f(g(x))].
\end{eqnarray*}
Using \eqref{3.18} and \eqref{3.19}, for any $t>0$ sufficiently small and $h'$ close to $h$, one has
\begin{eqnarray*}
&&\left\|\frac{1}{t}[f(g(x)+\triangledown g(x)(th')+o(t))-f(g(x)+\triangledown g(x)(th))]\right\|\\
&\leq&\frac{L}{t}\|t\triangledown g(x)(h'-h)+o(t)\|\rightarrow 0 \ (t\rightarrow 0^+, h'\rightarrow h).
\end{eqnarray*}
Then
\begin{eqnarray*}
\liminf_{t\rightarrow 0^+, h'\rightarrow h}\frac{1}{t}[\varphi(x,th')-\varphi(x)]&=&\liminf_{t\rightarrow 0^+}\frac{1}{t}[f(g(x)+t\triangledown g(x)(h))-f(g(x))]\\
&=&\mathbf{d}^+f(g(x), \triangledown g(x)(h)).
\end{eqnarray*}
This means that \eqref{3.17} holds. The proof is complete.  \hfill$\Box$\ \\[1pt]

\noindent{\it Proof of  \cref{theorem:th4}.} 
Let $\varphi:=f\circ g$. By the metric regularity of $g$ around $\bar x$, there exist $\kappa,\delta_0>0$ such that
\begin{equation}\label{4.17a}
  \mathbf{d}(x, g^{-1}(y))\leq \kappa\|g(x)-y\|\ \ \forall (x,y)\in \mathbf{B}(\bar x,\delta_0)\times \mathbf{B}(g(\bar x),\delta_0).
\end{equation}
Note that ${\rm bd}(S_f)\subseteq f^{-1}(0)$ and then
\begin{equation*}\label{4.27}
  {\rm bd}(\mathcal{S})\subseteq g^{-1}({\rm bd}(S_f))\subseteq g^{-1}(f^{-1}(0))=\varphi^{-1}(0).
\end{equation*}
By virtue of  \cref{lemma:lem0}
and \cref{proposition:pro4},
there is $\delta\in (0,\delta_0)$ sufficiently small such that $g^{-1}(S_f)$ has the Shapiro first order contact property on $\mathbf{B}(\bar x,\delta)$ and
\begin{equation}\label{4.29}
  \mathbf{T}^\mathbf{B}(g^{-1}(S_f), x)=\triangledown g(x)^{-1}(\mathbf{T}^\mathbf{B}(S_f), g(x))\ \ \forall x\in \mathbf{B}(\bar x,\delta)\cap g^{-1}(S_f).
\end{equation}
This means that $\mathcal{S}$ has the Shapiro first order contact property  around $\bar x$ as $\mathcal{S}=g^{-1}(S_f)$.

Define $\varPsi :X\times \mathbb{R}\rightarrow Y\times\mathbb{R}$ as
$$
\varPsi (x,r):=(g(x),r)\ \ \forall (x,r)\in X\times \mathbb{R}.
$$
Then one can verify that
$$
{\rm epi}(\varphi)=\varPsi ^{-1}({\rm epi}(f)).
$$
 For any $(x,y)\in \mathbf{B}(\bar x,\delta)\times\mathbf{B}(g(\bar x),\delta)$ and any $r,s\in(\varphi(\bar x)-\delta,\varphi(\bar x)+\delta)$, one has
$$
\varPsi ^{-1}(y,s)=g^{-1}(y)\times\{s\}
$$
and it follows from \eqref{4.17a} that
\begin{eqnarray*}
\mathbf{d}((x,r), \varPsi ^{-1}(y,s))&=&\mathbf{d}((x,r), g^{-1}(y)\times\{s\})\\
&=& \mathbf{d}(x, g^{-1}(y))+|r-s|\\
&\leq&\kappa\|g(x)-y\|+|r-s|\\
&\leq&(\kappa+1)\|\varPsi (x,r)-(y,s)\|,
\end{eqnarray*}
which implies that $\varPsi$ is metrically regular around $(\bar x,\varphi(\bar x))$. Thus, \cref{proposition:pro4} gives that $\varPsi ^{-1}({\rm epi}(f))$ has the Shapiro first order contact property at $(\bar x,\varphi(x))$ and so $\varphi$ has the epigraphical Shapiro first order contact property at $\bar x$ (thanks to ${\rm epi}(\varphi)=\varPsi ^{-1}({\rm epi}(f))$).

Finally, by applying  \cref{theorem:th2} and \cref{lemma:lem10},
one can obtain the proof of  \cref{theorem:th4}.
The proof is complete.\hfill$\Box$\\[1pt]


The following theorem follows immediately from \cref{theorem:th4}.
\begin{theorem}\label{th5}
Suppose that $f$ is continuous and $\bar x\in\mathcal{S}$ is such that $g$ is metrically regular around $\bar x$. Then composite-convex inequality \eqref{3.16} has a local error bound at $\bar x$ if and only if
\begin{empheq}[box =\mybluebox]{equation*}
\limsup_{x\stackrel{{\rm bd}(\mathcal{S})}\longrightarrow \bar x}\mathbf{e}\big(\{h\in X: \mathbf{d}^+f(g(x), \triangledown g(x)(h))\leq 1\},  \triangledown g(x)^{-1}(\mathbf{T}^\mathbf{B}(S_{f}, g(x))\big)<+\infty.
\end{empheq}
Further, one has the following accurate estimate for the local error bound modulus:
\begin{empheq}[box =\mybluebox]{equation*}
 \tau(\mathcal{S}, \bar x)=\limsup_{x\stackrel{{\rm bd}(\mathcal{S})}\longrightarrow \bar x}\mathbf{e}\big(\{h\in X: \mathbf{d}^+f(g(x), \triangledown g(x)(h))\leq 1\},  \triangledown g(x)^{-1}(\mathbf{T}^\mathbf{B}(S_{f}, g(x))\big).
\end{empheq}
\end{theorem}

The following example is to show the application of Theorem 4.5 when one verifies the error bound modulus of the composite-convex inequality.\\

\noindent{\bf Example 4.1.} Let $f(y):=y-1$ for all $y\in\mathbb{R}$, $g(x):=x^3$ for all $x\in\mathbb{R}$ and $\bar x:=1$. We consider the composite-convex inequality:
$$
(f\circ g)(x)\leq 0.
$$
Then one can verify that $\mathcal{S}=S_f=(-\infty, 1]$, $\mathbf{T}^\mathbf{B}(S_{f}, g(\bar x))=(-\infty,0]$ and
$$
\{h: \mathbf{d}^+f(g(\bar x), \triangledown g(\bar x)(h))\leq 1\}=(-\infty, \frac{1}{3}].
$$
This implies that
$$
\limsup_{x\stackrel{{\rm bd}(\mathcal{S})}\longrightarrow \bar x}\mathbf{e}\big(\{h\in X: \mathbf{d}^+f(g(x), \triangledown g(x)(h))\leq 1\},  \triangledown g(x)^{-1}(\mathbf{T}^\mathbf{B}(S_{f}, g(x))\big)=\frac{1}{3}.
$$
Hence Theorem 4.5 implies that $\mathcal{S}$ has a local error bound at $\bar x$ and moreover $ \tau(\mathcal{S}, \bar x)=\frac{1}{3}$.\hfill$\Box$\\

Finally, we give an example to show the local error bound may not be satisfied if the metric regularity assumption is dropped in Theorem 4.5.\\

\noindent{\bf Example 4.2.} Let $f(y):=y$ for all $y\in\mathbb{R}$, $g(x):=x^3$ for all $x\in\mathbb{R}$ and $\bar x:=0$. We consider the composite-convex inequality:
$$
(f\circ g)(x)\leq 0.
$$
Then $\mathcal{S}=S_f=(-\infty, 0]$, and one can verify that $g$ is not metrically regular at $\bar x$ since $\triangledown g(\bar x)=0$ is not surjective. However, for any $x>0$ sufficiently small, one has
$$
\frac{d(x, \mathcal{S})}{f(g(x))}=\frac{x}{x^3}\rightarrow +\infty, \ \ {\rm as} \ x\rightarrow 0^+.
$$
This means that $\mathcal{S}$ has no local error bound at $\bar x$, and thus the conclusions in Theorem 4.5. do not hold.\hfill$\Box$\\

\noindent{\bf Acknowledgements.} The authors are very grateful to the anonymous reviewers for their suggestions and comments that improved the presentation of this paper.

\section{Conclusions}

This paper is to establish primal characterizations of error bounds for a composite-convex inequality (that is defined by a composition of a convex function with a continuously differentiable mapping). To this aim, several primal necessary and/or sufficient conditions for a general inequality, under the assumption of Shapiro contact property, are given in terms of Bouligand tangent cones, directional derivatives and the Hausdorff-Pompeiu excess. Then it is proved that the composite-convex inequality satisfies the Shapiro contact property at the given point where the continuously differentiable mapping is metrically regular and thus primal characterizations of error bounds can be obtained. Our works actually extend the existing primal results on error bounds for the convex inequality to the non-convex case. The future work would be to investigate a broader class of non-convex inequality for which these primal results on error bounds are valid.

\end{document}